\documentclass[11pt]{article}

\usepackage{amsmath,amsthm,amsfonts}

\usepackage{fullpage}

\theoremstyle{plain}
\newtheorem{teo}{Theorem}[section]
\newtheorem{lem}[teo]{Lemma}
\newtheorem{coro}[teo]{Corollary}
\newtheorem{prop}[teo]{Proposition}
\theoremstyle{definition}
\newtheorem{ex}[teo]{Example}
\newtheorem{rem}[teo]{Remark}

\def\ra{\overline}

\def\CC{\mathbb{C}}
\def\O{\mathcal{O}}
\def\U{\mathcal{U}}
\def\Z{\mathcal{Z}}
\def\GL{\mathrm{GL}}
\def\SL{\mathrm{SL}}
\def\Sp{\mathrm{Sp}}

\def\ad{\mathrm{ad}}

\def\Hom{\mathrm{Hom}}

\def\Diff{\mathrm{Diff}}

\def\ZZ{\mathbb{Z}}

\def\wt{\widetilde}
\def\ot{\otimes}

\def\Ker{\mathrm{Ker}}

\def\Ext{\mathrm{Ext}}
\def\Tor{\mathrm{Tor}}
\def\Aut{\mathrm{Aut}}
\def\Pic{\mathrm{Pic}}
\def\Der{\mathrm{Der}}
\def\InnDer{\mathrm{InnDer}}
\def\VdB{\mathrm{V\!dB}}
\def\VDB{\mathrm{V\!dB}}
\def\op{\mathrm{op}}
\def\ch{\mathrm{ch}}

\def\cl#1{{\langle #1\rangle}}

\title{Hochschild duality, localization and smash products}
\author{Marco Farinati}

\begin{document}

\maketitle

\begin{abstract}
In this work we study the class of algebras satisfying a 
duality property with respect to Hochschild homology and 
cohomology, as in  \cite{VdB}.
 More precisely, we consider the class of algebras $A$ such that there exists
an invertible bimodule $U$ and an integer number $d$
with the property $H^{\bullet}(A,M)\cong
H_{d-\bullet}(A,U\ot_AM)$, for all $A$-bimodules $M$.
We will show that this class is closed under localization 
and  under smash products with respect to Hopf algebras satisfying also
the duality property.

We also illustrate the subtlety on dualities
with smash products developing in detail the example
 $S(V)\#G$, the crossed 
product of the symmetric algebra on a vector space and a 
finite group
acting linearly on $V$. 
\end{abstract}

\section*{Introduction}

The aim of this work is to study the class of algebras satisfying a 
duality property with respect to Hochschild homology and 
cohomology, as in  \cite{VdB}.
 More precisely, we consider the class of algebras $A$ such that there exists
an invertible bimodule $U$ and an integer number $d$
with the property $H^{\bullet}(A,M)\cong
H_{d-\bullet}(A,U\ot_AM)$, for all $A$-bimodules $M$.
We will show that this class is closed under localization 
(theorem \ref{teoloc}) and  under smash products 
(theorem \ref{teogal}). By localization we mean an algebra
 morphism $A\to B$ with the following two
properties: $B\ot_AB\cong B$ as $B$-bimodule,
and $B\ot_A-\ot_AB$ is exact.
For smash product, the philosophy is the following:
take $A$ an algebra in this class with dualizing bimodule $U$, 
and $H$ a Hopf algebra with dualizing bimodule $H$,
then $A\#H$ has dualizing bimodule $U\#H$ (see remark \ref{rem:smash}
for the definition of $U\#H$).

There is  a subtlety on dualities
with smash products, so the last section is devoted to develop
the simplest example
illustrating this:  the algebra $S(V)\#G$, the crossed 
product of the symmetric algebra on a vector space and a 
finite group
acting linearly on $V$. Given al algebra $A$ with dualizing module
$U_A\cong A$ and a Hopf algebra with dualizing bimodule 
isomorphic to $H$,
theorem \ref{teogal} says that $A\#H$ has a dualizing bimodule
isomorphic to $U_A\#H$. The subtlety is that, eventhow the
bimodule $U_A\cong A$ as $A$-bimodule, it may happens that
 $U_A\not\cong A$ as $H$-module, and so $U_A\#H\not\cong A\#H$ as
$A\#H$-bimodule. In the example of $S(V)$ and $G\subset \GL(V)$,
we show that the condition for $U_{S(V)}\cong S(V)$ as $G$-modules
is that $G\subset \SL(V)$, and consequently, homology
and cohomology will differ. In order to illustrate the duality, we compute
the cohomology of this example in two different ways.

The example of section 3 was motivated by a question
of Paul Smith, whether the methods used in \cite{AFLS} would apply
to $S(V)\#G$. The answer to that question is yes, and this calculation
has also motivated section 2. I am grateful to Jacques Alev to have transmitted
this question to me. I also want to thank Mariano Su\'arez \'Alvarez
for careful reading of this manuscript.

\subsection*{General notations}

Fix a field $k$ of characteristic zero, unadorned $\ot$ and $\Hom$ will
denote $\ot_k$ and $\Hom_k$.
If  $X$ is a graded vector space and  $n\in\ZZ$, we will denote
 $X[n]$ the same vector space but with its degree shifted by $n$.
For example, if $X$ is non-zero only in degree zero,
then  $X[n]$ is non-zero only in degree $n$.

For any $k$-algebra $B$ and $k$-symmetric bimodule $M$, 
the Hochschild homology and cohomology of $B$ with coeficients
in $M$ are $\Tor_{\bullet}^{B^e}(B,M)$ and 
$\Ext^{\bullet}_{B^e}(B,M)$, respectively, where
$B^e=B\ot B^{\op}$; they are denoted $H_{\bullet}(B,M)$ and 
$H^{\bullet}(B,M)$.
In the special case where $M=B$, we will also write
$HH_{\bullet}(B):=H_{\bullet}(B,B)$ and
$HH^{\bullet}(B):=H^{\bullet}(B,B)$.

The word ``module'' will mean ``left module''. All modules
will be $k$-symmetric, so that  $B$-bimodules is the same as
$B^e$-modules. A $B$-bimodule $P$ is called {\bf invertible} if there exists
another bimodule $Q$ such that $P\ot_BQ\cong B$
and $Q\ot_BP\cong B$. The set of isomorphism classes
of invertible $B$-bimodules
which are $k$-symmetric is denoted by $\Pic_k(B)$.

Finally, in section 3 there is some abuse of notation with  the symbol
$\det$. Some times it denotes the usual determinant function, and
some other times it denotes the 1-dimensional representation of
$\GL(V)$, or its restriction to some $G\subset\GL(V)$.
The meaning will be clear from the context.

\section*{The duality theorem of Van den Berg}

In  \cite{VdB}, the author proves a theorem
relating the Hochschild homology and cohomology of a certain
class of algebras.  
We will state this theorem in a way convenient for
our purposes:

\begin{teo}
\label{teo:dual}
(Theorem 3 of \cite{VdB}). Let $A$ be a $k$-algebra
which admits a finitely generated  projective $A^e$-resolution (for 
instance, this is the case if $A^e$ is noetherian) . The
following conditions  are equivalent:
\begin{enumerate}
\item There exists an invertible $A$-bimodule $U_A$, 
and an integer $d$ such that
$H^{\bullet}(A,M)\cong H_{d-\bullet}(A,U_A\ot_AM)$ for all $A^e$-modules $M$.
\item The projective dimension of $A$ as $A^e$-module is finite, and
$\Ext_{A^e}^n(A,A^e)=0$ for all $n\geq 0$ except for $n=d$ where
$U_A:=\Ext_{A^e}^n(A,A^e)$ is an invertible $A^e$-module.
\end{enumerate}
\end{teo} 

\section{Localization}

The general framework of this section
is the following:  $A\to B$ is a $k$-algebra map
such that
\begin{itemize}
\item The multiplication map induces an isomorphism of $B^e$-modules
$B\ot_AB\cong B$.
\item The functors $B\ot_A-$ and $-\ot_AB$ are exact.
\end{itemize}

We look for conditions on $B$ which, together
with the assumption that 
$A$ satisfies Van den Bergh's theorem, 
allow us to conclude that so does $B$.

\begin{lem} 
\label{lemapic}
Let $U\in \Pic(A)$ and $A\to B$ be such that $B\ot_AB\cong B$.
If $U\ot_AB\cong B\ot_AU$ as $A^e$-modules, then 
\begin{itemize}
\item $B\ot_AU \cong B\ot_AU\ot_AB$ as $B\ot A^{\op}$-modules;
\item $U\ot_AB \cong B\ot_AU\ot_AB$ as $A\ot B^{\op}$-modules; and
\item  $B\ot_AU$ is a $B^e$-module in a natural way,
$B\ot_AU \in \Pic(B)$, its inverse is $B\ot_AU^{-1}\ot_AB$,
and $U^{-1}\ot_AB\cong B\ot_AU^{-1}$ as $A^e$-modules.
\end{itemize}
\end{lem}
\begin{proof}
The first isomorphism is the composition:
\[B\ot_A(U\ot_AB)\cong B\ot_A(B\ot_AU)=
(B\ot_AB)\ot_AU\cong B\ot_AU\]
The second one is similar.

Now let $U^{-1}$ be the inverse of $U$ in $\Pic(A)$,
so that  $U\ot_AU^{-1}\cong
U^{-1}\ot_AU\cong A$. Let 
us see  that $B\ot_AU^{-1}\ot_AB$ is the  inverse of
$B\ot_AU$:
\[
\begin{array}{rl}
(B\ot_AU)\ot_B(B\ot_AU^{-1}\ot_AB)&\cong
(U\ot_AB)\ot_BB\ot_AU^{-1}\ot_AB\\
&\cong U\ot_AB\ot_AU^{-1}\ot_AB\\
&\cong B\ot_AU\ot_AU^{-1}\ot_AB\\
&\cong B\ot_AA\ot_AB\\
&\cong B\ot_AB\\
& \cong B,
\end{array}
\]
and
\[
\begin{array}{rl}
(B\ot_AU^{-1}\ot_AB)\ot_B(B\ot_AU)&\cong
B\ot_AU^{-1}\ot_AB\ot_AU\\
&\cong B\ot_AU^{-1}\ot_AU\ot_AB\\
&\cong B\ot_AA\ot_AB\\
&\cong B\ot_AB\\
& \cong B.
\end{array}
\]
\end{proof}
A bimodule $U$ such that there is an isomorphism $B\ot_AU\cong
U\ot_AB$ of $A^e$-modules will be said to {\bf commute with $B$}.

\begin{ex}
Let $g\in\Aut_k(A)$ be such that it admits an extension
$\wt{g}\in\Aut_k(B)$, i.e. $\wt{g}(a)=g(a)$ for all $a\in A$.
Then the element $Ag\in\Pic(A)$ commutes with $B$.
In particular, $U=A$ commutes with $B$.
\end{ex}
\begin{proof}
Let $g$ be such an element and consider $Ag\in\Pic(A)$.
There is an isomorphism of $B\ot A^{\op}$-modules
\[
\begin{array}{rcl}
B\ot_AAg&\to& B\wt{g}\\
b\ot ag&\mapsto& ba\wt{g}
\end{array}
\]
On the other hand, one can define an isomorphism of $A\ot B^{\op}$-modules
\[
\begin{array}{rcl}
Ag\ot_AB&\to& B\wt{g}\\
ag\ot a\wt{g}&\mapsto& a\wt{g}(b)\wt{g}
\end{array}
\]
In particular, $Ag\ot_AB$ and $B\ot_AAg$ are isomorphic as $A^e$-modules.
\end{proof}
\begin{ex}
Let $g\in\Aut_k(A)$ be such that there exists no element
$\wt{g}\in\Aut_k(B)$ extending it. Then the bimodule
$Ag$ doesn't commutes with $B$.
\end{ex}
\begin{proof}
Assume 
$B\ot_AAg\cong Ag\ot_AB$ as $A^e$-modules. From lemma
\ref{lemapic} it follows that $B\ot_AAg\in\Pic(B)$. But, as a left
$B$-module, 
$B\ot_AAg\cong B$, and it is well-known that if an element
$U\in\Pic(B)$ is such that ${}_BU\cong {}_BB$, then it
 is  of the form
$B\alpha$ for some $\alpha\in\Aut_k(B)$, the automorphism $\alpha$ 
being defined
up to inner automorphism. In particular, for $a\in A$ one
has that
$
g(a)=u\alpha(a)u^{-1}$
for some $u\in\U(B)$. Denoting
$\wt{g}:=u\alpha(-)u^{-1}$ we see that we have found an automorphism
extending $g$, thus a contradiction.
\end{proof}
\begin{rem}
Let $A\to B$ be such that
$B\ot_AB\cong B$.
If $M$ is a left $B$-module, then $M\cong B\ot_AM$ as a left $B$-module.
If  $N$ is another left $B$-module, then
$\Hom_B(M,N)=\Hom_A(M,N)$.
\end{rem}

\begin{proof}
Using the hypothesis on $B$, we see that
\[M
\cong B\ot_BM
\cong (B\ot_AB)\ot_BM
\cong B\ot_A(B\ot_BM)
\cong B\ot_AM;\]
it follows then that
\[
\Hom_B(M,N)
\cong\Hom_B(B\ot_AM,N)
\cong\Hom_A(M,N).
\]
\end{proof}

\begin{teo}
\label{teoloc}
Let $A\in VdB(d)$ with dualizing bimodule $U$, 
and
 $A\to B$ be a morphism of $k$-algebras such that
\begin{enumerate}
\item the functors $B\ot_A-$ and $-\ot_AB$ are exact; 
\item the canonical map induced by multiplication
$B\ot_AB\to B$ is an isomorphism; and
\item $B\ot_AU\cong U\ot_AB$ as $A^e$-modules.
\end{enumerate}
Then 
$B\in VdB(d)$ with dualizing bimodule $B\ot_AU\cong B\ot_AU\ot_AB$.
\end{teo}
Notice that if $U=A$, then  condition 3 is automatically
satisfied, and the
dualizing bimodule associated to $B$ is $B$.

\begin{proof}
By theorem \ref{teo:dual}, it is enough to show that the projective dimension
of $B$ as $B^e$-module is finite, that $B$ admits a resolution by means of finitely generated $B^e$-projectives.
 and that  $\Ext_{B^e}^n(B,B^e)=B\ot_AU\ot_AB$ and it vanishes
elsewhere.

Let $P_{\bullet}$ be a finite resolution of $A$ as $A^e$-modules, with
$P_n$ projective and finitely generated as $A^e$-modules. Since
$B\ot_A-$ and $-\ot_AB$ are exact, the complex
$B\ot_AP_{\bullet}\ot_AB$ is a resolution of
$B\ot_AA\ot_AB\cong B$,
and so $B$ also has a finite resolution. 
The bimodules $B\ot_AP_n\ot_AB$ are clearly $B^e$-finitely generated
and projective.

In order to compute $\Ext^{\bullet}_{B^e}(B,B^e)$ one can 
use this particular resolution, and consequently
\[
\Ext^{\bullet}_{B^e}(B,B^e)= 
H^{\bullet}(\Hom_{B^e}(B\ot_AP_{\bullet}\ot_AB,B^e))
\cong  H^{\bullet}(\Hom_{A^e}(P_{\bullet},B^e))\]
We claim that if $P$ is $A^e$-projective finitely generated, then
\[
\Hom_{A^e}(P_{\bullet},B^e)\cong
B\ot_A\Hom_{A^e}(P_{\bullet},A^e)\ot_AB
\] 
For that, consider the class of $A^e$-modules $P$ such that
$\Hom_{A^e}(P_{\bullet},B^e)\cong
B\ot_A\Hom_{A^e}(P_{\bullet},A^e)\ot_AB$. This class is closed
under direct summands and finite sums, so it is enough to show 
our claim that
the module $A^e$ is in it, and that is clear.
Using this isomorphism one gets
\[  H^{\bullet}(\Hom_{A^e}(P_{\bullet},B^e))
\cong  H^{\bullet}(B\ot_A\Hom_{A^e}(P_{\bullet},B^e)\ot_AB) \]
and by flatness this is the same as
$B\ot_A  H^{\bullet}(\Hom_{A^e}(P_{\bullet},B^e))\ot_AB 
=B\ot_AU[d]\ot_AB$.
\end{proof}

\begin{ex}\label{a1}
We can take $A=A_1(k)=k\{x,y\}/\cl{[x,y]=1}$, 
$B=k\{x,x^{-1},y\}/\cl{[x,y]=1}$.
This example is a particular case of the following:
\end{ex}
\begin{ex}{\bf Normal localization:}
Let $A$ be an algebra and $x\in A$ such that the set
$\{1,x,x^2,x^3,\dots\}$ satisfies the Ore
conditions. Take $B=A[x^{-1}]$.
If $M$ is a right $A$-module, then as $k[x]$ modules we have an isomorphism
$M\ot_AB\cong M\ot_{k[x]}k[x,x^{-1}]$. This shows that $A\to B$ is flat.
It is also clear that $B\ot_AB\cong B$, in the same way as
$k[x^{\pm1}]\ot_{ k[x]} k[x^{\pm1}]\cong k[x^{\pm1}]$.
\end{ex}
\begin{ex}
Another  generalization of example \ref{a1} is the following situation:
let $\O(X)$ be the algebra of functions on an affine variety
$X$, and let $U$ be an affine open subset of $X$. Let $A=\Diff(X)$ be
the algebra of algebraic differential operators on $X$ and similarly
$B=\Diff(U)$. Since $B=\O(U)\ot_{\O(X)}\Diff(X)$, the map $A\to B$ is flat,
and  $B\ot_AB=B$.
If $A$ satisfies the theorem of Van den Bergh, then so it does $B$.
\end{ex}

Next section, we will study the behavior of the duality property
with respect to smash products.

\section{Smash products}

In this section $H$ is a hopf algebra such that $H\in VdB(d)$ with
dualizing bimodule $U_H=H$, 
$A\in VdB(d')$ is an $H$-module algebra with dualizing bimodule $U_A$, 
and $B:=A\#H$.
We will prove (see theorem \ref{teogal})
that $B\in VdB(d+d')$, with dualizing bimodule $U_B=U_A\#H$ (see remark
\ref{rem:smash} for the definition of $U\#H$).

\begin{lem}
If $H$ is a Hopf algebra, then $H\in VdB(d)$ with dualizing bimodule $H$ if and only if
$\Ext^{\bullet}_H(k,M)\cong \Tor_{\bullet-d}(k,M)$ for all left $H$-modules $M$.
\end{lem}
\begin{proof}
Let $M$ be a left $H$-module, then $M_{\epsilon}$ is the $H^e$-module with
right
action defined by $m.h:=\epsilon(h)m$ for all $m\in M$ and $h\in H$. If $H\in VdB(d)$,
it follows that
\[
\Ext^{\bullet}_H(k,M)=
H^{\bullet}(H,M_{\epsilon})\cong
H_{d-\bullet}(H,M_{\epsilon})
= \Tor_{\bullet-d}(k,M)\]
On the other direction, if $X$ is an $H^e$-module, then $X^{\ad}$ is the same underlying vector
space but with left $H$ action defined by
$h\cdot_{\ad}x:=h_1xS(h_2)$. With this structure (see for instance 
\cite{St}) one has
\[
H^{\bullet}(H,X)=
\Ext^{\bullet}_H(k,X^{\ad})\cong
 \Tor_{\bullet-d}(k, X^{\ad})\cong
H_{d-\bullet}(H,M_{\epsilon})
= \Tor_{\bullet-d}(k,M)\]
\end{proof}

\begin{ex}
Let $G$ be a finite group such that $\frac{1}{|G|}\in k$. The
Reynolds operator $e=\frac{1}{|G|}\sum_{g\in G}g$ induces an isomorphism
$M_G\cong M^G$
for any $G$-module $M$.
This implies that $k[G]\in VdB(0)$ with $U_{k[G]}=k[G]$.
This example can be easily  generalized in the following direction: 
\end{ex}
\begin{ex}
Let $H$ be a semisimple unimodular Hopf algebra, so that
$H$ admits a {\em central} integral $e\in H$ 
satisfying
\[he=\epsilon(h)e,\ \qquad \epsilon(e)=1. \]
Then $H\in\VdB(0)$ with $U_H=H$.
It is known (see Radford, \cite{R} theorem 4)
that the Drinfel'd double of a finite dimensional hopf algebra is unimodular.
If $K$ is a finite dimensional Hopf algebra and $D(K)$ is the 
Drinfel'd double, again by a result of Radford (\cite{R} proposition 7)
$D(K)$ is semisimple if and only if $K$ is 
semisimple and cosemisimple. Taking $K=k[G]$ where $G$ is a 
non-commutative group with $|G|^{-1}\in k$, we get $H:=D(K)$ a non 
commutative not cocommutative semisimple unimodular Hopf algebra.
\end{ex}
\begin{proof}
Let $H$ be a unimodular semisimple Hopf algebra, and let $e\in H$
be as above. We will show that $\Hom_H(k,M)\cong k\ot_HM$.
If $M$ is a left $H$-module, then
\[
\Hom_H(k,M)\cong\{m\in M\ / \ hm=\epsilon(h)m\}=:M^H.\]
 It is clear that
every element of the form $em$ belongs to $M^H$ because 
\[
h(em)=(he)m=\epsilon(h)em;
\]
but
 if $m\in M^H$, then 
\[em=\epsilon(e)m=m,
\] 
so $M^H$ coincides with the image of the 
multiplication by $e$. Let us consider the map
\[
\begin{array}{rcl}
e:M&\to& M^H\\
m&\mapsto & em.
\end{array}
\]
The elements of the form $hm-\epsilon(h)m$ belong to the kernel of this map, so it factors through
$M_H:=M/\cl{hm-\epsilon(h)m}$. Now the map $M^H\to M_H$ defined by
$m\mapsto \ra{m}$ defines an inverse, because in $M_H$, every element
$m=\epsilon(e)m$ is equivalent to $em$. We have shown that $H\in VdB(0)$.
\end{proof}

\begin{ex}\label{ejemplokx}
The algebra $H=k[x]$ is a Hopf algebra with $\Delta(x)=x\ot 1+1\ot x$.
It belongs to the Class $VdB(d)$ with $U_H=H$.
\end{ex}
\begin{proof}
Write $k[x]^e=k[x]\ot k+\cong =k[x,y]$, and consider the Koszul resolution
\[0\to k[x,y] \to 
 k[x,y] \to k[x]\to 0\]
where the first map is the multiplication by $(x-y)$ and the
second map is the evaluation $x=y$.
Applying the functor $\Hom_{k[x,y]}(-,k[x,y])$ on obtain the complex
\end{proof}
\[ 0\to \Hom_{k[x,y]}(k[x,y],k[x,y]) \to  \Hom_{k[x,y]}(k[x,y],k[x,y]) \to  0\]
where the map is again multiplication by $x-y$. This complex identifies with 
\[ 0\to k[x,y] \to  k[x,y] \to  0\]
but notice that now the grading increases to the right, so the 
homology is $k[x,y]/(x-y)\cong k[x]$ in degree one, zero
elsewhere, and we conclude that $k[x]\in \VdB(1)$.

\begin{ex} 
The algebra $k[x]$ admits a finitely generated $k[x]^e$-projective resolution;
this fact implies a K\"unneth formula for Hochschild cohomology, and
so the algebra $k[x_1,\dots,x_n]\in VdB(n)$, with
$U_{k[x_1,\dots,x_n]}=
k[x_1,\dots,x_n]$.
\end{ex}

\begin{ex}
The Hopf algebra $k[x_1^{\pm 1},\dots,x_d^{\pm1}]=k[\ZZ^n]$, 
belongs to the class $\VdB(d)$, because as an algebra, it is a localization of
$k[x_1,\dots,x_d]$. Also
\[U_ {k[x_1^{\pm 1},\dots,x_d^{\pm1}]}
=U_ {k[x_1,\dots,x_d}]\ot_{k[x_1,\dots,x_d]}
k[x_1^{\pm 1},\dots,x_d^{\pm1}]=
k[x_1^{\pm 1},\dots,x_d^{\pm1}]\]
\end{ex}

\begin{rem}
\label{rem:smash}
Let $A$ be an $H$-module algebra and $U\in \Pic_k(A)$ such that
$U$ is also an $H$-module, with the compatibility property
\[
h(aub)=h_1(a)h_2(u)h_3(b)
\]
for all $a,b \in A$, $h \in H$, and $u \in U$.
Let $U^{-1}:=\Hom_A(U,A)$; this is also an $H$-module satisfying the same
compatibility condition.
If $U\#H$ is the abelian group $U\ot H$ with $A\#H$-bimodule
structure
given by
\[
(a\#h)(u\ot k):=
(ah_1(u)\ot h_2k)
\]
\[
(u\ot k)(a\#h)=
(uk_1(a)\ot k_2h),
\]
then $U\#H\in\Pic_k(A\#H)$, and  its inverse is $U^{-1}\#H$.
If $M$ is left $A\#H$-module, then
\[
(U\#H)\ot_{A\#H}M\cong U\ot_AM
\] 
as $A\#H$-modules, where the 
$A\# H$-module structure on $U\ot_AM$ is the one induced by the obvious left $A$-structure
and the diagonal $H$-structure.
\end{rem}
\begin{proof}
We will only  exhibit an isomorphism $U\#H\ot_{A\#H}U^{-1}\#H\to A\#H$.
Let us denote  by $\cl{\ ,\ }$ the evaluation map $U\ot_AU^{-1}\to A$; notice
that $\cl{\ ,\ }$ is $H$-linear. For
$u\in U$, $v\in U^{-1}$, $h$ and $k\in H$, define
\[
\begin{array}{rcl}
U\#H\ot_{A\#H}U^{-1}\#H &\to &A\#H\\
(u\ot h)\ot (v\ot k)&\mapsto & \cl{u,h_1(v)}h_2k.
\end{array}
\]
\end{proof}

\begin{teo}
\label{teogal}Let $H\in VdB(d)$ be a Hopf algebra with $U_H=H$.
If  $A$ is an $H$-module algebra with $A\in VdB(d)$, then
$A\# H \in VdB(d+d')$ with $U_{A\# H}=U_{A}\#H$.
\end{teo}
\begin{proof}
Let $B$ be $A\# H$.
In \cite{St}, the author shows that, for a $B$-bimodule $M$, there
is a spectral sequence converging to $H^{\bullet}(B,M)$ whose second term
is $\Ext^p(k,H^q(A,M))$. Similarly, there is a spectral sequence with $E^2$ term
equal to $\Tor_p^H(k,H_q(A,M))$ converging to $H_{\bullet}(B,M)$.

Now consider $M=B^e$, and let us compute $H^{\bullet}(B,B^e)$. 
First, one notes the following isomorphism of left
$A^e$-modules:
\[B^e\cong A^e \ot V,\]
where $V$ is the vector space $H\ot H$.

Using Stefan's spectral sequence, one has
\[
E_2^{pq}=\Ext_H^p(k,H^q(A,B^e))
=\Ext_H^p(k,H^q(A,A^e\ot V)).
\]
Since $A\in VdB(d')$, it follows that
\[\begin{array}{rcl}
H^{\bullet}(A,A^e\ot V))
&\cong& H_{d'-\bullet}(A,U\ot_AA^e\ot V)\\
&\cong& H_{d'-\bullet}(A,U\ot_AA^e)\ot V\\
&\cong& H^{\bullet}(A,A^e)\ot V\\
&\cong& U[d]\ot V.
\end{array}\]
This implies first that the spectral sequences degenerates at this 
step, and
consequently, there is an isomorphism
\[
H^{\bullet}(B,B^e)\cong \Ext_H^{*-d'}(k,U\ot V)
\]
Recall that $V=H\ot H^{op}$; we have to consider it as 
$H$-module with 
the adjoint action. Now we use the fact that $H\in VdB(d)$, with
$U_H=H$, so
$H^{\bullet}(H,X)\cong H_{d-\bullet}(H,X)$ for all $H$-bimodules $X$. In particular,
for a left $H$-module $X$, one can consider the bimodule
$X_{\epsilon}$, and this gives the formula
\[
\Ext_H^{\bullet}(k,X)
=H^{\bullet}(H,X_{\epsilon})
\cong H_{d-\bullet}(H,X_{\epsilon})
=\Tor^H_{d-\bullet}(k,X).\]
This formula implies that
\[
H^{\bullet}(B,B^e)\cong \Ext_H^{*-d'}(k,U_A\ot V)\cong
\Tor^H_{d'+d-\bullet}(k,U_A\ot V).\]

On the other hand,
$H_{\bullet}(B,U_A\ot_AB^e))=
H_{\bullet}(B,(U_A\#H)\ot_BB^e))$
can be computed using a spectral sequence
whose second term is
\[
\begin{array}{rcl}
\Tor_{\bullet}^H(k,H_{\bullet}(A,U_A\ot_AB^e))
&=&
\Tor_{\bullet}^H(k,H_{\bullet}(A,U\ot_A(A^e\ot V))\\
&=&\Tor_{\bullet}^H(k,U_A\ot V).
\end{array}
\]
This spectral sequence collapes giving an isomorphism
\[
H_{\bullet}(B,U_A\ot_AB^e)
\cong \Tor_{\bullet}^H(k,U_A\ot V)
\]
In particular, 
\[
H^{\bullet}(B,B^e)\cong H_{d+d'-\bullet}(B,U\ot_AB^e)
\]
and 
\[
\begin{array}{rcl}
H^{d+d'}(B,B^e)
&=& H_{0}(B,U\ot_AB^e)\\
&=&H_{0}(B,(U\#H)\ot_BB^e)\\
&=&H_{0}(B,(U\#H)\ot B)\\
&=&U\#H.
\end{array}
\]
\end{proof}

\begin{coro}
With the notations of the above theorem, assume $U=A$ as $A$-bimodules
and $H$-modules, then
\[
H^{\bullet}(B,M)\cong H_{d+d'-\bullet}(B,M)
\]
for all $A\#H$-bimodules $M$.
\end{coro}

\begin{ex}Let $A\in VdB(d)$,  $D\in\Der_k(A)$, and
write the Ore extension $B=A[t,D]$. 
This algebra $B$  coincides with 
$A\#k[t]$ where the $k[t]$-module action on $A$
is given by $t.a=D(a)$, $B\in VdB(d+1)$. For $A=k[x]$ and
$D=\frac{\partial}{\partial x}$ one
obtains the known result that $A_1(k)\in VdB(2)$.
\end{ex}

\begin{ex} Let $0\neq q\in k$, then 
$B=k\{x^{\pm 1},y^{\pm 1}\}/\cl{yx=qxy}\in VdB(2)$. Indeed, this algebra is isomorphic to
$k[x^{\pm 1}]\# k[y^{\pm 1}]$ where the $H$-module structure
on $k[x^{\pm 1}]$ is given by $y.x=qx$.
\end{ex}

\begin{ex}
Let $A$ be an algebra and $G$ a finite group of automorphism of
$A$. If $A\in VdB(d)$, then $A\#G\in \VdB(d)$.
\end{ex}

\noindent {\bf Warning:} It can happen that $A$ is such that $U_A\cong A$
as $A$-bimodule, but $U_A\not\cong A$ as $H$-module.
It is easy to show an example of this situation when $H=k[G]$.

One can first observe the following caracterization of the
$A^e\#G$-structures on a $A$-bimodule isomorphic to $A$:
\begin{prop}
\label{propstr}
Let $U$ be an $A^e$-bimodule isomorphic to $A$. 
The set of all possible $A^e\# G$-module structures on
$U$, modulo $A^e\#G$-isomorphism,
is parametrized by $H^1(G,\U\Z(A))$, the first cohomology of $G$
with coeficients in the (multiplicative) abelian group of 
units of the center of $A$.
\end{prop}
\begin{proof}
Fix an isomorphism $A\cong U$ and let $u$ be the image of $1$
in $U$. Hence $U=Au=uA$, and  moreover, $au=ua$ for
all $a\in A$.
One has to define a $G$-action on $U$ such that,
for all $a,b\in A$ and $v\in U$, the following identity holds
\[g(avb)=g(a)g(v)g(b).\]
Since the bimodule $U$ is generated by $u$, it is clear that it is 
only necesary to define $g(u)$. 
The element $g(u)$ must belong to $U$, so it is of the form $a_gu$ for some $a_g$ in $A$. 
But 
\[au=ua\] for all $a\in A$, and applying $g$ one obtains
\[
ag(u)=g(u)a,\ \ \forall a\in A,\]
and so 
\[aa_gu=a_gua=a_gau.\] 
It follows that $a_g$ must belong
to the center of $A$. Also, every element of $U$ is of the form
\[
ag(u)=aa_gu,\]
 so $a_g$ must be a unit. 
We have then shown that the assignment $g\mapsto a_g$
must be a map from $G$ into $\U(\Z(A))$.

If one wants associativity, the identity 
\[
g(h(u))=(gh)(u),\   \forall g,h\in G.
\]
is required, so
\[\begin{array}{rcl}
g(h(u))&=&g(a_hu)\\
&=&g(a_h)a_gu\\
&=&(gh)(u)\\
&=&a_{gh}u.
\end{array}\]
But $u$ is a basis of $U$ with 
respect to the left $A$-structure,
so
\[g(a_h)a_g=a_{gh}. \]

On the other hand, it is clear that  an assignment $g\mapsto a_g$
from $G$ into the units of center of $A$ satisfying the above cocycle condition
defines a $G$-action compatible with the $A$-bimodule structure.

Now assume that $U$ has two $G$-actions that are isomorphic.
Let us denote them by
$g._1(u)=a_gu$, and $g._2(u)=b_gu$, and call $U_1$ and
$U_2$ the bimodule $U$ with the first and the second $G$-structure, respectively.

If $\phi:U_1\to U_2$ is an isomorphism of $A^e\#G$-modules,
then
the image of $u$ is some element $\lambda u$, where
$\lambda\in A$. Moreover, $\lambda$ is a unit 
because $\phi$ is an isomorphism, and $\lambda\in
\Z(A)$ because $\phi$ is $A^e$-linear.

Now  $G$-linearity
means that
\[
\begin{array}{rcl}
\phi(g._1u)&=&\phi(a_gu)\\
&=&\lambda a_gu,
\end{array}
\]
but also 
\[
\begin{array}{rcl}
\phi(g._1u)
&=&g._2\phi(u)\\
&=&g._2(\lambda u)\\
&=&g(\lambda)g._2 u\\
&=&g(\lambda)b_g u,
\end{array}
\]
so we deduce
\[b_g=\lambda g(\lambda^{-1})a_g,\]
and the two assignments differ by a coboundary.
 
\end{proof}

Despite proposition \ref{propstr}, for an algebra $A\in VdB$,
the dualizing  bimodule $U$ is a very particular one, namely
$U_A=Ext_{A^e}^d(A,A^e)$. The following is an example showing
(without calculating $H^1(G,\U\Z(A))$ that $U$ is isomorphic to
$A$ as $A^e$ bimodule, but not as $G$-module:

\begin{ex}\label{ejemplodet}
Let $V$ be a finite dimensional vector space, $A=S(V)$, and
$G\subset \GL(V)$ a finite group.
We claim that 
\[\Ext_{A^e}^{\bullet}(A,A^e)=A\ot\det{}^{-1}[d],
\] 
where $d=\dim(V)$,
and $\det{}^{-1}$ is the dual of the  determinant representation $\Lambda^d V$. 
Namely, 
$\det{}^{-1}$ is a one dimensional $k$-vector space, 
if $w\in\det{}^{-1}$ is a nonzero element,  $g\in G$, and $a\in A$, then
 the $G$-action is given by
\[g(a\ot w)=g(a)\det(g|_V)^{-1}\ot w.\]
We conclude that $U_A\cong A$ as $A^e\#G$-modules if and only if $G\subset \SL(V)$.
\end{ex}

\begin{proof}
Let $g\in G$, and choose a basis $\{x_1,\dots, x_d\}$ of $V$ which diagonalizes
$g$. Notice that $S(V)=\ot_{i=1}k[x_i]$, and this tensor product
is $g$-equivariant with the diagonal action. The
 K\"unneth formula is $g$-equivariant, so we only need to  prove 
the following lemma:

\begin{lem}\label{lemadet}
If $A=k[x]$ and $g$ is the automorphism of $A$ determined by
$g(x)=\lambda x$, then
$\Ext_{A^e}^{\bullet}(A,A^e)=A[1]$, and 
the action of $g$ is given by multiplication
by $\lambda^{-1}$.
\end{lem}
\noindent{\it Proof of the lemma:} 
It was shown in example \ref{ejemplokx} that 
$k[x]\in\VDB(1)$, let us compute the $g$-action 
on 
\[
H^1(k[x],k[x,y])=\Der(k[x],k[x,y])/\InnDer(k[x],k[x,y]).
\]
If $D:k[x]\to k[x,y]$ is a derivation, then $D$ is determined by
its value $D(x)$ on $x$, and this gives the isomorphism
\[
\begin{array}{rclr}
\Der(k[x],k[x,y])&\cong& k[x,y]  \\
D&\mapsto &D(x).&
\end{array}
(\dag)
\]

If $p\in k[x,y]$, the inner derivation
$[p,-]$ takes in $x$ the value
\[
\begin{array}{rcl} 
[p,x]&=&p(x,y)y-xp(x,y)\\
&=&(x-y)p(x,y).
\end{array}
\]
This shows that, under the isomorphism (\dag),  $\InnDer\cong (x-y)k[x,y]$, obtaining
\[
\begin{array}{rcl}
H^1(A,A^e)
&=&\Der(A,A^e)/\InnDer(A,A^e)\\
&\cong&\dfrac{k[x,y]}{(x-y)k[x,y]}\\
&\cong & k[x].
\end{array}
\]
In order to compute the action of $g$ on $H^1$ we recall that,
if $D$ is a derivation, then
\[
g.D=g(D(g^{-1}(-)),
\]
 so
\[
\begin{array}{rcl}
(g.D)(x)&=& g(D(g^{-1}x))\\
&=&g(D(\lambda^{-1} x))\\
&=&\lambda^{-1}g(D(x)),
\end{array}
\]
and if $D(x)\in k$ (this is always the case modulo an inner derivation)
we get
\[
(g.D)(x)=\lambda^{-1} D(x).
\]
\end{proof}
Back to the example $A=S(V)$ and $G\subset\GL(V)$ a finite subgroup, we see 
that $S(V)\#G\in \VDB(\dim(V))$ but $U_{S(V)\#G}\cong S(V)\#G$ if
and only if $G\subset\SL(V)$.
This example shows a situation where $H^{\bullet}(B,M)=H_ {d-\bullet}(B,U\ot_BM)$ with 
$U\neq B$. In particular, $H^{\bullet}(B)\cong H_{\bullet}(B,U)$,
which  needs not be equal to $H_{d-\bullet}(B)$, and in fact
it is different.

\section{The example $S(V)\#G$}

We finish with  a computation of the  homology and cohomology of 
$S(V)\#G$.

Let $k$ be a field, $V$ a finite dimensional $k$-vector space,
$G$ a finite subgroup of $\GL(V,k)$, $A=S(V)$, and we will asume that 
$\frac1{|G|}\in k$. 
For simplicity we will also asume that $k$ has a primitive 
$|G|$-th root of 1. This condition is not really necessary
because of the following reason: consider $\xi$ a primitive
$|G|$-root of unity in the algebraic closure of $k$ and
let $K$ be $k(\xi)$ the field generated by $k$ and $\xi$. 
One can view $G$ inside $\GL(V\ot K,K)$, and
consider it acting on $A\ot K=S_K(V\ot K)$. A descend property
of the Hochschild homology and cohomology  with respect to this
change of the base field
assures that the dimension over $K$ of the (co)homology
of the extended algebra is the same as the dimension over $k$ 
 of the (co)homology of the original one.

If $g\in G$, $V^g=\{x\in V \ /\ g(x)=x\}$. As $g$-module,
$V^g$ admits a unique complement in $V$, we will call it
$V_g$. We have  $V=V^g\oplus V_g$ as $g$-modules, and this decomposition 
is canonical.

\subsection{Homology of $S(V)\#G$}

\begin{teo}
\label{teo1}
With the notations as in the above paragraph,
denote $\cl{G}$ the set of conjugacy classes of $G$, and
for $g\in G$ let $\Z_g$ be the centralizer of $g$ in $G$,
so that  $\Z_g=\{h\in G\ /\ hg=gh\}$. The Hochschild homology
of $S(V)\#G$ is given by:
\[
H_n(S(V)\#G)=
H_n(S(V),S(V)\#G)^G=
\bigoplus_{\cl{g}\in\cl{G} }(S(V^g)\ot \Lambda^n (V^g))^{\Z_g}\]
where $\Lambda^n(V^g)$ is the homogeneous component of degree $n$ 
of the exterior algebra on $V^g$.
\end{teo}

\begin{proof}

With the hypothesis on the characteristic and the order of the group,
the spectral sequence of \cite{St} gives the following
isomorphism:
\[
\begin{array}{rcl}
H_n(S(V)\#G)&=&H_n(S(V),S(V)\#G)^G\\
&=&\bigoplus_{\cl{g}\in \cl{G} }H_n(S(V),S(V)g)^{\Z_g},
\end{array}
\]
valid for any $k$-algebra of the type $A\#G$.
Since $V=V^g\oplus V_g$, it follows that 
\[
S(V)\cong S(V^g)\ot S(V_ g)\]
as algebras, and
\[
S(V)g\cong S(V^g)\ot S(V_ g)g\]
as  $S(V)$-bimodules. Using 
the  K\"unneth formula one gets
\[
H_n(S(V),S(V)g)^{\Z_g}
=\bigoplus_{p+q=n}
\left( H_p(S(V^g))\ot H_q(S(V_g),S(V_g)g) \right) ^{\Z_g}\]
By the 
Hochschild-Kostant-Rosenberg theorem, or directly by computing
using a Koszul type resolution, one see that,
if $W$ is a finite dimensional $k$-vector space,
\[
H_n(S(W))=\Omega^n(S(W))=S(W)\ot \Lambda^n W.\]

The  homology with coeficients is computed in the following lemma:
\begin{lem}
\label{lema:coef}
$H_{\bullet}(S(V_g),S(V_g)g)=k[0]$
with  trivial $\Z_g$-action.
\end{lem}
\begin{proof} 
Let  $h\in \Z_g$. One can diagonalize simultaneously 
$h$ and $g$ in $V_g$. If $\{x_1,\dots x_k\}$ is a basis of eigenvectors
of both $h$ and $g$, then the algebra $S(V_g)$ is isomorphic
to
\[k[x_1,\dots x_k]=\bigotimes_{i=1}^kk[x_i]
\] and
\[S(V_g)g=k[x_1,\dots x_k]g=\bigotimes_{i=1}^kk[x_i]g_i,
\] where
$g_i$ acts on $x_i$ by multiplication of the corresponding
eigenvalue og $g$.
Notice also that $h$ acts on each $x_i$ by multiplication by some
$\lambda_i'$, because
$x_i$ is also an eigenvector of $h$.

Using the  K\"unneth formula again, one gets:
\[
H_{\bullet}(S(V_g),S(V_g)g)=\bigotimes_i H_{\bullet}(k[x_i],k[x_i]g_i).
\]
Let us now make the explicit computation for the algebra
$k[x]$, $g$ acting by
$x\mapsto \lambda x$, and $h$ acting by $x\mapsto \lambda'x$. 

Consider, as in example \ref{ejemplokx}, the resolution
of $k[x]$ as $k[x]$-bimodule
\[0\to k[x]\ot k[x]\to k[x]\ot k[x]\to k[x]\to 0.\]
Here the first morphism is given by
$p\ot q\mapsto px\ot q-p\ot xq$ and the second one is the multiplication
map.

By tensoring with  $k[x]g$ over $k[x]^e$, one gets the
complex
\[0\to k[x]g\to k[x]g\to 0\]
with differential
\[pg\mapsto pgx-xpg=px(\lambda-1)g,\]
whose homology is $H_{\bullet}(k[x].k[x]g)$.
The fact that  $\lambda\neq 1$ implies that the differential is injective
and the image equals $xk[x]g$, so  $H_1=0$ and $H_0=k$. It is clear that
$h$ acts trivially on $H_0$, and the proof of the lemma is
complete.
\end{proof}

The sum 
\[
H_n(S(V),S(V)g)^{\Z_g}
=\bigoplus_{p+q=n}
\left(H_p(S(V^g))\ot H_q(S(V_g),S(V_g)g) \right) ^{\Z_g}\]
reduces to
\[
H_n(S(V),S(V))^{\Z_g}=
(S(V^g)\ot\Lambda^n(V^g) ) ^{\Z_g}\]
and the proof of the theorem is finished.
\end{proof}

\begin{ex}
Let $k=\CC$, $V=\CC^2$, $G$ a finite subgroup of
$\SL(2,\CC)$. Then
\[
\begin{array}{rcl}
H_0(S(V)\#G)&=&
S(V)^G\oplus \CC^{\#\{\cl{g}\neq 1\}}\\
H_1(S(V)\#G)&=&
(S(V)\ot V)^G\\
H_2(S(V)\#G)&=&
(S(V)\ot\Lambda^2(V))^G= S(V)^G\\
H_n(S(V)\#G)&=&0\ \forall n>2
\end{array}
\]
\end{ex}

\subsection{Cohomology: direct computation}

The formula
\[
H^n(S(V)\#G)=
H^n(S(V),S(V)\#G)^G=
\bigoplus_{\cl{g}\in \cl{G} }H^n(S(V),S(V)g)^{\Z_g}\]
is also valid. Using $S(V)=S(V^g)\ot S(V_ g)$,
and the K\"unneth formula one gets
\[\begin{array}{rcl}
H^n(S(V),S(V)g)^{\Z_g}
&=&\underset{p+q=n}{\bigoplus}
\left( H^p(S(V^g),S(V^g))\ot H^q(S(V_g),S(V_g)g) \right) ^{\Z_g}\\
&=&\underset{p+q=n}{\bigoplus}
\left( S(V^g)\ot \Lambda^p((V^g)^{\bullet})\ot H^q(S(V_g),S(V_g)g) \right) ^{\Z_g}.
\end{array}\]
Here we have used the isomorphism
\[
\begin{array}{rcl}
H^{\bullet}(S(W),S(W))
&=&\Lambda^{\bullet}_{S(W)}\Der(S(W))\\
&=&S(W)\ot\Lambda^{\bullet} W^*.
\end{array}
\]
Now we need the analogue of the lemma \ref{lema:coef} for cohomology,
whose proof is the same as lemma \ref{lemadet}.
\begin{lem}
Let $A=k[x]$, $g,h$ the automorphisms determined by
$g(x)=\lambda x$ and $h(x)=\mu x$, with $\lambda\neq 1$. Then
$H^{\bullet}(A,Ag)=k[1]$, and the action of $h$ is given by multiplication
by $\mu^{-1}$.
\end{lem}

\begin{coro}If we denote by $d_g=\dim_k(V_g)$ then
\[
H^{\bullet}(S(V_g),S(V_g)g)=\det|_{V_g}^{-1}[d_g].
\]
This is an isomorphism of  of $\Z_g$-modules.
\end{coro}
\begin{proof}
From the fact that $g$ and $h$ commute, one can choose a basis
$\{x_1,\dots,x_n\}$ of eigenvectors of both $g$ and $g$.
The corollary follows from the K\"unnet formula, and the
Lema above applied
to $S(V)=\otimes_{i=1}^nk[x_i]$.
\end{proof}

We have obtained the following formula:
\begin{teo}
\label{teodirecto}
\[
H^{\bullet}(S(V)\#G)=
\bigoplus_{\cl{g}\in\cl{G}}
\left(S(V^g)\ot\Lambda^{\bullet}((V^g)^{\bullet})\ot
\det|_{V_g}^{-1}[d_g]\right)^{\Z_g}.
\]
\end{teo}

\subsection{Cohomology: computation using duality}

Using theorem \ref{teogal} for $H=k[G]$ (see example \ref{ejemplodet}), 
we know that
\[
\begin{array}{rcl}
H^{\bullet}(A\#G)&=&H^{\bullet}(A\#G,(U_A\#G)\ot_{A\#G}A\#G)\\
&=&H_ {d-\bullet}(A\#G,U_A\#G)\\
&=&H_ {d-\bullet}(A\#G,(A\ot\det{}^{-1})\#G).
\end{array}
\]
Using Stefan's spectral, this is the same as
\[
\begin{array}{rcl}
H_ {d-\bullet}(A,(A\ot\det{}^{-1})\#G)^{G}
&=& \bigoplus_{\cl{g}\in\cl{G}}
H_ {d-\bullet}(A,(A\ot\det{}^{-1}(V)).g)^{\Z_g}\\
&=&\bigoplus_{\cl{g}\in\cl{G}}
(H_ {d-\bullet}(A,A.g)\ot\det{}^{-1}(V))^{\Z_g}.
\end{array}\]
Now the same techniques of writing $V=V^g\oplus V_g$ apply, and
we obtain
\[
\begin{array}{rcl}
\underset{\cl{g}\in\cl{G}}{\bigoplus}
(H_ {d-\bullet}(A,A.g)\ot\det{}^{-1})^{\Z_g}
&=&
\underset{\cl{g}\in\cl{G}}{\bigoplus}
(H_ {d-\bullet}(S(V^g))\ot\det{}^{-1})^{\Z_g}\\
&=&
\underset{\cl{g}\in\cl{G}}{\bigoplus}
(S(V^g)\ot\Lambda^{d-\bullet}(V^g)\ot\det{}^{-1})^{\Z_g}.
\end{array}\]
The difference between this formula and that of Theorem
\ref{teodirecto}, having $\det$ or $\det|_{V_g}$ is explained by the fact
that in
\ref{teodirecto}, one has
also $\Lambda^{\bullet}((V^g)^{*})$, while here one has
$\Lambda^{d-\bullet}(V^g)$. The multiplication map induces a morphism
of $\Z_g$-modules
\[\Lambda^{\bullet}(V^g)\ot\Lambda^{\dim(V^g)-\bullet}(V^g) 
\to \Lambda^{\dim(V^g)}V^g=\det|_{V^g},\]
and as a consequence one has an isomorphism of $\Z_g$-modules
\[
\Lambda^{\bullet}(V^g)^*\cong \Lambda^{\dim(V^g)-\bullet}(V^g) 
\ot\det|_{V^g}^{-1}.\]
So we get the same after noticing that $\det=\det|_{V^g}\ot\det|_{V_g}$.

\begin{ex}
Let $k=\CC$, $V=\CC^2$, $G$ a finite subgroup of $\SL(2,\CC)$.  In this
case, homology and cohomology is the same:
\[
\begin{array}{rcl}
H^0(S(V)\#G)&=&
S(V)^G\\
H^1(S(V)\#G)&=&
(S(V)\ot V)^G\\
H^2(S(V)\#G)&=&
S(V)^G\oplus \CC^{\#\{\cl{g}\neq 1\}}\\
H^n(S(V)\#G)&=&0\ \forall n>2
\end{array}
\]
\end{ex}

\begin{ex}
Let $G=C_2=\{1,t\}$ the cyclic group of order two. Let $k$ be a field
of $\ch(k)\neq 2$, $A=k[x]$ with $t$ acting on $A$ by
$x\mapsto -x$.
Using theorem \ref{teo1} one gets

\[
\begin{array}{rcccl}
H_0(A\#G)&=&
A^G\oplus k&=&k[x^2]\oplus k\\
H_1(A\#G)&=&
(A\ot k.dx)^G&=&k[x^2]xdx \\
H_n(A\#G)&=&0&&\forall n>1
\end{array}
\]
On the other hand,
\[
\begin{array}{rcccl}
H^0(A\#G)&=&
A^G&=&k[x^2]\\
H^1(A\#G)&=&
(A\ot k.\partial_x)^G \oplus 
\left(
\frac{\Der(A,At)}{\InnDer(A,At)}\right)^{C_2}
&=&k[x^2]x\partial_x\oplus 0 \\
H_n(A\#G)&=&0&&\forall n>1
\end{array}
\]
 In this example, homology and cohomology are not the same.
The cohomology is $k[x^2]$-free, while the homology has torsion.
\end{ex}

In the above example, we see that the cohomology is a ``part'' of the homology.
The same phenomenon happens in the following:
\begin{ex} Let $W=k^{n}$, consider $S_n$ acting on $W$ by permutation of 
the coordinates,
and let \[
V=\{(1,1,\dots,1)\}^{\perp}:=\{(x_1,\dots,x_n)\in W\ / \ 
\sum_{i=1}^nx_i=0\}.\]
We claim that 
\[
H^{\bullet}(S(V)\#S_n)=
H^{\bullet}(S(V), S(V)\#A_n)^{S_n},
\]
where $A_n$ denote as usual the subgroup
of even permutations.
\end{ex}

In fact, we can prove an analogous formula in the following general setting:
\begin{ex}
Let $G\subset\GL(V)$ be a finite subgroup,  $S:=G\cap\SL(V)=
\Ker(\det:G\to k^{\times})$, and $C:=\det(G)\subset k^{\times}$.
Then
\[
H_{\bullet}(S(V)\#G)=
\bigoplus_{w\in C}
\left(
\bigoplus_{\cl{g}\in\cl{G},\ det(g)=w}
H_{\bullet}(S(V),S(V)g)^{\Z_g}
\right),\]
and each of this summands is non zero, while in cohomology, there
are only the terms corresponding to $w=1$:
\[
H^{\bullet}(S(V)\#G)=
\bigoplus_{\cl{g}\in\cl{G},\ \det(g)=1}H^{\bullet}(S(V),S(V)g)^{\Z_g}.\]
In particular
\[
H^{\bullet}(S(V)\#G)=
H^{\bullet}(S(V), S(V)\#S)^G \]
and
\[ H^{\bullet}(S(V)\#G)\neq 
 H_{d-\bullet}(S(V)\#G).\]
\end{ex}
\begin{proof}
The formula for the homology is just noticing that the set $\cl{G}$
can be split into smaller pieces, parametrized by the values 
of the determinant.
To see that each summand is non-zero we make them explicit. Using
theorem \ref{teo1} we know that:
\[
H_{\bullet}(S(V),S(V)g)^{\Z_g}
=(S(V^g)\ot\Lambda^{\bullet}V^g)^{\Z_g}.
\]
Even if $V^g=0$, one always has the element $1\in
(S(V^g)\ot\Lambda^{\bullet}V^g)^{\Z_g}$.

The interesting part is the formula for the cohomology.
Recall from the  duality formula that
\[H^{\bullet}(S(V),S(V)g)\cong 
\det{}^{-1} \otimes H_{d-\bullet}(S(V),S(V)g).\]
If one shows that 
$H_{\bullet}(S(V),S(V)g)$ is a trivial $g$-module, then, for
$\det(g)\neq 1$ we will have
\[\begin{array}{rcl}
\left(\det{}^{-1} \otimes H_{d-\bullet}(S(V),S(V)g)\right)^{\Z_g}
&\subseteq&
\left(\det{}^{-1} \otimes H_{d-\bullet}(S(V),S(V)g)\right)^{g}\\
&=&
(\det{}^{-1})^g \otimes H_{d-\bullet}(S(V),S(V)g)\\
&=&0.
\end{array}\]
So let us see that $H_{\bullet}(S(V),S(V)g)$ has trivial $g$-action. 
For that, write $V=V^g\oplus V_g$, then 
$H_{\bullet}(S(V),S(V)g) \cong
H_{\bullet}(S(V^g))\otimes H_{\bullet}(S(V),S(V)g))$. 
Clearly
 $H_{\bullet}(S(V^g))$ is a trivial $g$-module, and 
 $H_{\bullet}(S(V),S(V)g))$ also has trivial $g$-action in virtue of
lema \ref{lema:coef}.
\end{proof}
\begin{rem}The equality between homology and cohomology depends not 
only on $G$, but on the representation. For example, given an arbitarry
 finite subgroup $G\subset\GL(V)$, we can consider the action 
on $V$ and on $V^*$, and $G$ will act
symplectically
on $W=V\oplus V^*$. In this case we have 
\[
G\hookrightarrow \Sp(W)\subset
\SL(W),\]
 so that
\[
H^{\bullet}(S(W)\#G)=H_{\dim(W)-\bullet}(S(W)\#G).
\]
\end{rem}

\end{document}